\documentclass[12pt]{amsart}
\usepackage{amsmath,amsthm,amsfonts,amssymb}
\begin{document} 
\newcommand{\B}{{\mathbb B}}
\newcommand{\C}{{\mathbb C}}
\newcommand{\N}{{\mathbb N}}
\newcommand{\Q}{{\mathbb Q}}
\newcommand{\Z}{{\mathbb Z}}
\renewcommand{\P}{{\mathbb P}}
\newcommand{\F}{{\mathcal F}}
\newcommand{\R}{{\mathbb R}}
\newcommand{\rc}{\subset}
\newcommand{\rank}{\mathop{rank}}
\newcommand{\trace}{\mathop{tr}}
\newcommand{\dimc}{\mathop{dim}_{\C}}
\newcommand{\Lie}{\mathop{Lie}}
\newcommand{\Auto}{\mathop{{\rm Aut}_{\mathcal O}}}
\newcommand{\alg}[1]{{\mathbf #1}}
\newtheorem*{definition}{Definition}
\newtheorem*{claim}{Claim}
\newtheorem{corollary}{Corollary}
\newtheorem*{Conjecture}{Conjecture}
\newtheorem*{SpecAss}{Special Assumptions}
\newtheorem{example}{Example}
\newtheorem*{remark}{Remark}
\newtheorem*{observation}{Observation}
\newtheorem*{fact}{Fact}
\newtheorem*{remarks}{Remarks}
\newtheorem{lemma}{Lemma}
\newtheorem*{proposition}{Proposition}
\newtheorem*{theorem}{Theorem}
\title{%
Surface foliations with compact complex leaves are holomorphic
}
\author {J\"org Winkelmann}
\begin{abstract}
Let $X$ be a compact complex surface with a real foliation.
If all leaves are compact complex curves, the foliation
must be holomorphic.

Soit $X$ une surface complexe compacte munie d'un feuilletage
reel. Supposons que tous les feuilles sont de courbes complexes
compactes. Alors le feuilletage est holomorphe.
\end{abstract}
%
\address{%
J\"org Winkelmann \\
 Institut Elie Cartan (Math\'ematiques)\\
 Universit\'e Henri Poincar\'e Nancy 1\\
 B.P. 239\\
 F-54506 Vand\oe uvre-les-Nancy Cedex\\
 France
}
\email{jwinkel@member.ams.org\newline\indent{\itshape Webpage: }%
http://www.math.unibas.ch/\~{ }winkel/
}
\maketitle
\section{The result}
\begin{theorem}
Let $X$ be a compact complex smooth surface, and $\F$ a real foliation
on $X$ such that all the leaves are compact complex curves.

Then $\F$ is in fact a holomorphic foliation.
\end{theorem}
\begin{proof}
We use the theory of cycle space resp.~Chow schemes.%
\footnote{One might also use the theory of Hilbert schemes.
Basically, Chow and Hilbert schemes differ only in their
treatment of non-reduced subvarieties and this is of no
relevance to our problem.}
By this theory there exists a ``cycle space'' $C$ with countably
many irreducible components $C_i$ all of which are compact,
and a universal space $U\subset C\times X$ such that for every
compact complex irreducible subspace $Z$ of $X$ there is a point
$t\in C$ such that $U\cap\left(\{t\}\times X\right)=\{t\}\times Z$.

This is classical for projective varieties, but actually works for
arbitrary compact complex manifolds
by work of Barlet \cite{B}. 
We will need the compactness
of the components $C_i$. This is true for all surfaces as well as for
all K\"ahler manifolds. However, there are higher-dimensional non-K\"ahler
manifolds for which $C_i$ need not be compact
(see e.g.~\cite{W}, cor.~4.11.3).

Because there are only countably many connected components $C_i$ there
must be a component $C_0$ such that the collection of subvarieties
parametrized by this $C_0$ includes uncountably many of the leaves
$L_\alpha$ of the foliation $\F$. Let $U_0$ denote the corresponding
component of $U$, and let $p:U_0\to C_0$ and $\pi:U_0\to X$ denote
the projections.
Because everything is compact, $\pi(U_0)$ is an irreducible closed
analytic subset of $X$. Since it contains infinitely many distinct
curves, we must have $\pi(U_0)=X$.

Now fix a point $q\in C_0$ such that the parametrized subvariety
$Z_q=\pi(p^{-1}(q))$ is one of the leaves of $\F$.

Let $L$ be an arbitrary leaf of $\F$ other than $Z_q$. Then 
$L\cap Z=\emptyset$
and therefore $L\cdot Z=0$ in the sense of intersection theory.
Now for every $s\in C_0$ the corresponding curve $Z_s$ is numerically
equivalent to $Z_q$. Hence $L\cdot Z_s=0$ for all $s\in C_0$
(and every leaf $L$ of $\F$).
On the other hand, there must be a parameter $s\in C_0$ such that
$L\cap Z_s\ne\emptyset$, because the union of the $Z_s$ covers $X$
(i.e., $\pi(U_0)=X$).
As a consequence, every leaf $L$ must coincide with an irreducible
component of $Z_s$ for some $s\in C_0$.

As a next step we want to show that actually all the $Z_s$ are
irreducible. There are several ways to do this.

One way is to use the
fact that all $Z_s$ are connected, and therefore
their irreducible components intersect while the leaves are by
definition disjoint.

Here we use a different method based on discussing volumes.
We observe that for every K\"ahler form $\omega$ the volume
$\int_{Z_s}\omega$ is locally constant as a function of $s\in C_0$.
This implies that for almost all leaves $L$ we have $\int_L\omega=K=
\int_{Z_q}\omega$. Now for the foliation $\F$ the volume of the leaves
$\int_L\omega$ must vary continuously in dependence of the leaf $L$.
Hence it can not jump and therefore it
equals $K$ for every leaf $L$. However, if some $Z_s$ is reducible,
an irreducible component of $Z_s$ evidently has volume smaller than $K$.
Therefore such a component  can not occur as leaf of $\F$. 
It follows that for a reducible
$Z_s$ the equality $Z_s\cdot L=0$ implies that $Z_s$ has empty
intersection with {\em every} leaf $L$. Since the union of all leaves
covers $X$, this is absurd. Hence there is no reducible $Z_s$.

Thus we have seen: Every $Z_s$ (with $s\in C_0$) is an irreducible 
curve on $X$ and every leaf $L$ of $\F$ equals one of the curves $Z_s$.

Conversely, assume that $s\in C_0$ and consider the curve $Z_s$:
Since the leaves of $\F$ cover all of $X$, there is a leaf $L$ with
$Z_s\cap L\ne\emptyset$. With $Z_s\cdot L=0$ it follows that $Z_s=L$.

We have thus established that there is a one-to-one correspondance between
the leaves of $\F$ and the curves $C_s$ parametrized by $s\in C_0$.

It follows that $\F$ is a holomorphic foliation.
\end{proof}
\section{Higher dimensions}
The result is not valid in higher dimensions.
Indeed, let $M$ be a real 4-manifold with an anti-self dual Riemannian
metric $g$, e.g. a compact real $4$-dimensional torus with flat
metric.
Then we have a ``twistor space'' $X$ with a projection $\pi:X\to M$
such that $X$ is a compact complex three-dimensional manifold
and all the fibers of $\pi$ are compact complex curves of genus $0$
whose normal bundle is ${\mathcal O}(1)^{\oplus 2}$ (\cite{AHS}).
Because the normal bundle of the fibers is
not holomorphically trivial, it is clear that there is no way to define a
complex structure on $M$ for which is $\pi$ is holomorphic.
It follows in this way also that these curves are not leaves of
a holomorphic foliation on $X$.

\section{Non-compacxt case}
Compactness of the surface is crucial,
as can be seen by the following example
which is due to J.J.~Loeb.

We consider the non-compact surface
$X=\P_2(\C)\setminus
\P_2(\R)$.
Each $x\in\P_2(\C)$ corresponds to a complex line in $\C^3$.
Each vector $v\in\C^3$ decomposes into a real part and an imaginary
part: $v=u+iw$, $u,w\in\R^3$. The real vectors $u,w$ are linearly
independent unless $x\in\P_2(\R)$. Thus we obtain a map from $X$
to $\P_2(\R)^*$ by mapping each $x=[u+iw]\in X$ to the
real hyperplane of $\R^3$ spanned by $u$ and $w$.
This yields a real-analytic map $F:X\to\P_2(\R)^*$ whose fibers
are complex curves: If $H=\left<u,w\right>_{\R}\in\P_2(\R)^*$, then
\[
F^{-1}(H)=\{ [v]\in\P_2(\C)\setminus\P_2(\R)
: v\in\left<u,w\right>_{\C} \}
\]
Let $v\in\P_2(\R)$. Then for every real hyperplane $H$ in $\R^3$ containing
$v$ we can find a sequence $x_n\in X$ with $F(x_n)=H$ and
$\lim x_n=v$. Therefore $F$ can nowhere extended to $\P_2(\R)$
as a continuous map.

Now assume that the foliation defined by $F$ is holomorphic.
Each $F$-fiber has two components (the fibers are isomorphic to
$\P_1(\C)\setminus\P_1(\R)$). Hence $F$ lifts to a map $\tilde F$
to $S^2$, the $2:1$-covering of $\P_2(\R)$.
If $F$ defines a holomorphic foliation, this map $\tilde F$
must be holomorphic for some complex structure on $S^2$.
But then $\tilde F$ would be a meromorphic function,
and meromorphic functions extend through totally real
submanifolds like $\P_2(\R)$ in $\P_2(\C)$ while we have
seen that $F$ and $\tilde F$ do not even extend as topological
maps.

Thus $F$ defines a foliation on $X$ whose leaves are all
isomorphic to $H^+=\{z\in\C:\Im(z)>0\}$,
but this foliation is not holomorphic.

\end{document}